# RIFFLE SHUFFLES OF DECKS WITH REPEATED CARDS

By Mark Conger and D. Viswanath[1]

*University of Michigan*

By a well-known result of Bayer and Diaconis, the maximum entropy model of the common riffle shuffle implies that the number of riffle shuffles necessary to mix a standard deck of 52 cards is either 7 or 11—with the former number applying when the metric used to define mixing is the total variation distance and the latter when it is the separation distance. This and other related results assume all 52 cards in the deck to be distinct and require all 52! permutations of the deck to be almost equally likely for the deck to be considered well mixed. In many instances, not all cards in the deck are distinct and only the sets of cards dealt out to players, and not the order in which they are dealt out to each player, needs to be random. We derive transition probabilities under riffle shuffles between decks with repeated cards to cover some instances of the type just described. We focus on decks with cards all of which are labeled either 1 or 2 and describe the consequences of having a symmetric starting deck of the form $1, \ldots, 1, 2, \ldots, 2$ or $1, 2, \ldots, 1, 2$. Finally, we consider mixing times for common card games.

**1. Introduction.** The connection between examples and concepts in probability theory is a particularly close one. That examples derived from the question "How many shuffles mix a deck of cards?" have featured prominently in the development of the convergence theory for Markov chains by Persi Diaconis and others can be seen in this light. This article deals with riffle shuffling, which is the most common way of shuffling cards.

There are $2^n$ ways to cut a deck of $n$ cards into two packets and then riffle them together since a card that ends up in the $i$th position can be dropped by either the left hand or the right hand. The maximum entropy model assigns equal probability to all these $2^n$ riffle shuffles. More generally, the maximum entropy model assigns equal probability to all $a^n$ $a$-shuffles, with

Received October 2003; revised January 2005.

[1]Supported in part by a fellowship from the Sloan Foundation.

an $a$-shuffle being a way to cut a deck into $a$ packets and then riffle them together. Several equivalent descriptions of the $a$-shuffle have been given by Bayer and Diaconis [2]. The $a$-shuffle with $a = 2$ is also described by Epstein [5] who calls it the amateur shuffle.

We will refer to elements of the group $S_n$ of permutations of $\{1, 2, \ldots, n\}$ as shuffles. If $\pi \in S_n$ and $\pi(i) = j$, then by convention the shuffle $\pi$ sends the card in the $i$th position to the $j$th position. The number of descents of $\pi$ is defined as the number of positions $1 \leq i < n$ at which $\pi(i) > \pi(i+1)$. Bayer and Diaconis [2] proved that the probability that an $a$-shuffle results in a shuffle $\pi$ with $d$ descents is given by $\frac{1}{a^n}\binom{a+n-d-1}{n}$. The picture below shows a 3-shuffle of 6 cards.

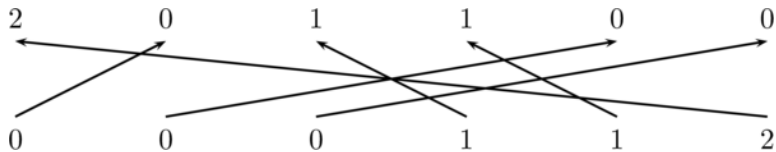

The bottom line indicates that the 0th, 1st and 2nd packets in the cut have 3, 2 and 1 cards, respectively. The top line indicates that the 1st, 2nd, 3rd, 4th, 5th and 6th cards in the shuffled deck are dropped from the 2nd, 0th, 1st, 1st, 0th and 0th packets, respectively. If the numbers are ignored, the arrows alone depict a shuffle. In a depiction of a shuffle such as the one above, a descent corresponds to a crossing between arrows that originate at adjacent positions. The shuffle depicted above has 2 descents, and therefore, according to Bayer and Diaconis [2], the probability that an $a$-shuffle results in the shuffle depicted above is $\frac{1}{a^6}\binom{a+3}{6}$.

In nearly all of the literature on card shuffling, the cards in a deck are assumed to be distinct. We allow cards to be indistinguishable. In our notation, both $1, 1, 2, 1$ and $1^2, 2, 1$ denote the deck with two cards labeled 1 above a card labeled 2 above a card labeled 1. Let $a_1, a_2, \ldots, a_n$ be a deck. When it is shuffled using $\pi \in S_n$, the deck obtained is $a_{\pi^{-1}(1)}, a_{\pi^{-1}(2)}, \ldots, a_{\pi^{-1}(n)}$. We define $\pi(D_1; D_2)$ as the set of shuffles $\pi \in S_n$ such that $\pi$ applied to $D_1$ results in $D_2$. The *descent polynomial* of the shuffles from the starting deck $D_1$ to the ending deck $D_2$ is defined as

$$\sum_{\pi \in \pi(D_1; D_2)} x^{\mathrm{des}(\pi)},$$

where $\mathrm{des}(\pi)$ is the number of descents in $\pi$. For example, the descent polynomial of shuffles from $1, 1, 2, 2$ to $1, 2, 2, 1$ is $2x + 2x^2$.

If the descent polynomial of the shuffles from a deck $D_1$ with $n$ cards to a deck $D_2$ is $\sum_{d=0}^{n-1} c_d x^d$, the probability that an $a$-shuffle of $D_1$ results in the

deck $D_2$ is

$$p_a = \sum_{d=0}^{n-1} \frac{c_d}{a^n} \binom{a+n-d-1}{n}, \tag{1.1}$$

a formula obtained by summing over all the shuffles in $\pi(D_1; D_2)$. The system of linear equations (1.1) can be inverted to obtain

$$c_d = p_{d+1}(d+1)^n - p_d d^n \binom{n+1}{1}$$
$$+ p_{d-1}(d-1)^n \binom{n+1}{2} - \cdots + (-1)^d p_1 1^n \binom{n+1}{d}, \tag{1.2}$$

for $0 \le d < n$. It is possible to pass back and forth between the transition probabilities $p_a$ and the descent polynomial of $\pi(D_1; D_2)$ using (1.1) and (1.2).

In Section 2 we deduce efficient recursions for the descent polynomial from the starting deck $D_1$ to the ending deck $D_2$ when either $D_1$ or $D_2$ is a sorted deck of the form $1^{n_1}, 2^{n_2}, \ldots, k^{n_k}$. We derive a formula for the transition probabilities when $D_1 = 1, 2, \ldots, k, x^n$ without using the descent polynomial. Sections 3 and 4 consider the starting decks $D_1 = 1^n, 2^n$ and $D_1 = (1,2)^n$. Section 5 summarizes mixing times for card games obtained using results in the preceding sections and Monte Carlo simulations.

Although decks with repeated cards do not seem to have been considered, the work of Diaconis, McGrath and Pitman [4], Fulman [7] and Lalley [11] on cycle decompositions, and of Fulman [8] on increasing subsequences are in a somewhat similar vein. The thesis of Reyes [13] has new results, as well as many references related to other types of shuffles.

**2. Transition probabilities.** We begin with a recursive algorithm to obtain the descent polynomial of shuffles from $1^{n_1}, 2^{n_2}, \ldots, h^{n_h}$ to a deck $D$ which has the same $n_1 + n_2 + \cdots + n_h$ cards but in a different order. Each of the numbers $n_1, n_2, \ldots, n_h$ is a positive integer. The transition probability under an $a$-shuffle can be obtained using the descent polynomial and (1.1).

We assume the starting deck to be $1^{n_1}, 2^{n_2}, \ldots, h^{n_h}$, which is in sorted order. We denote by $D(i,c)$ the position of the $i$th card labeled $c$ in $D$. For example, if $D = 1, 2, 1, 1, 2, 2$, then $D(2, 1) = 3$. The deck obtained from $D$ by keeping the cards labeled $1, 2, \ldots, e$ and by discarding cards with all other labels will be denoted $D_{1e}$. Similarly, the deck obtained from $D$ by keeping the cards with labels $f, \ldots, h$ and discarding other cards will be denoted by $D_{fh}$. We assume $1 \le \cdots \le e < f \le \cdots \le h$ and that there is no card whose label is in-between $e$ and $f$, or, equivalently, $f = e + 1$. Let $N = n_1 + n_2 + \cdots + n_h$, $N_1 = n_1 + n_2 + \cdots + n_e$ and $N_2 = N - N_1$. Then $N$, $N_1$ and $N_2$ equal the number of cards in $D$, $D_{1e}$ and $D_{fh}$, respectively.

Consider the set of all shuffles $\pi$ from the sorted deck $1^{n_1}, 2^{n_2}, \ldots, h^{n_h}$ to $D$ such that $\pi(1) = D(i,1)$ and $\pi(N) = D(j,h)$, where $1 \leq i \leq n_1$ and $1 \leq j \leq n_h$. The number of these shuffles with $d$ descents is set equal to the coefficient of $x^d$ to define the polynomial $p_{i,j}(x)$.

To obtain a recursion for $p_{i,j}(x)$, consider the set of shuffles from the sorted deck $1^{n_1}, 2^{n_2}, \ldots, e^{n_e}$ to $D_{1e}$ and the set of shuffles from the sorted deck $f^{n_f}, \ldots, h^{n_h}$ to $D_{fh}$. Define $q_{i,k}(x)$, for $1 \leq i \leq n_1$ and $1 \leq k \leq n_e$, as the polynomial in which the coefficient of $x^d$ equals the number of shuffles $\pi$ with $d$ descents belonging to the first set which satisfy $\pi(1) = D_{1e}(i,1)$ and $\pi(N_1) = D_{1e}(k,e)$. The polynomial $r_{l,j}(x)$, for $1 \leq l \leq n_f$ and $1 \leq j \leq n_h$, is defined similarly, with the coefficient of $x^d$ equal to the number of shuffles $\pi$ with $d$ descents in the second set which satisfy $\pi(1) = D_{fh}(l,f)$ and $\pi(N_2) = D_{fh}(j,h)$. Then the following recursive relationship holds:

$$(2.1) \qquad p_{i,j}(x) = \sum_{k,l} q_{i,k}(x) r_{l,j}(x) x^{\varepsilon(k,l)}.$$

The indices $k$ and $l$ vary over $1 \leq k \leq n_e$ and $1 \leq l \leq n_f$. The exponent $\varepsilon(k,l)$ is 0 if $D(k,e) < D(l,f)$ and 1 if $D(k,e) > D(l,f)$.

To prove (2.1), we consider a bijection between $\pi(1^{n_1}, 2^{n_2}, \ldots, h^{n_h}; D)$ and $\pi(1^{n_1}, \ldots, e^{n_e}; D_{1e}) \times \pi(f^{n_f}, \ldots, h^{n_h}; D_{fh})$. Let the shuffle $\pi$ map to the pair of shuffles $(\pi_1, \pi_2)$ under this yet to be defined bijection. If position $i$ is occupied by a card labeled $\delta$ in the starting deck $1^{n_1}, 2^{n_2}, \ldots, h^{n_h}$ and $\pi(i) = D(j, \delta)$, then $\pi_1(i) = D_{1e}(j, \delta)$ if $1 \leq \delta \leq e$ and $\pi_2(i - N_1) = D_{fh}(j, \delta)$ if $f \leq \delta \leq h$, by definition of the bijection. To complete the proof of (2.1), we relate the number of descents of $\pi$ to the number of descents of $\pi_1$ and $\pi_2$. The number of descents of $\pi$ equals the sum of the number of descents of $\pi_1$ and $\pi_2$ if $\pi(N_1) = D(k,e)$, $\pi(N_1 + 1) = D(l,f)$ and $D(k,e) < D(l,f)$. However, if $D(k,e) > D(l,f)$, the sum must be incremented by 1.

The base case of the recurrence (2.1) occurs when the starting deck has cards of only one type. Take this deck to be $1^n$. The coefficient of $x^d$ in $p_{i,j}(x)$ is then equal to the number of shuffles $\pi \in S_n$ with $\pi(1) = i$, $\pi(n) = j$ and with $d$ descents. The number of shuffles $\pi \in S_n$ with $d$ descents is defined as the Eulerian number $\langle {n \atop d} \rangle$ [9]. Given a permutation of $\{1, 2, \ldots, n-1\}$ with $d$ or $d-1$ descents, the number $n$ can be inserted in $d+1$ or $n-d$ places, respectively, to obtain a permutation of $\{1, 2, \ldots, n\}$ with $d$ descents. Thus, as shown in [9], consideration of the insertion of the number $n$ into a permutation of the numbers $1, 2, \ldots, n-1$ gives the recurrence

$$(2.2) \quad \begin{aligned} \left\langle {n \atop d} \right\rangle &= (d+1) \left\langle {n-1 \atop d} \right\rangle + (n-d) \left\langle {n-1 \atop d-1} \right\rangle \qquad \text{if } n > 0, \\ \left\langle {0 \atop 0} \right\rangle &= 1, \qquad \left\langle {0 \atop d} \right\rangle = 0 \qquad \text{if } d \neq 0. \end{aligned}$$

The modified Eulerian number $\langle {n \atop d} \rangle_i$ is defined as the number of $\pi \in S_n$ with $\pi(1) = i$ and $d$ descents. If $d = 0$, $\langle {n \atop d} \rangle_i$ is 0 if $i > 1$ and 1 if $i = 1$. Consideration of the insertion of $n$ into a permutation of the numbers $1, 2, \ldots, n-1$ that begins with $i$ gives the recurrence

$$\left\langle {n \atop d} \right\rangle_i = (d+1) \left\langle {n-1 \atop d} \right\rangle_i + (n-d-1) \left\langle {n-1 \atop d-1} \right\rangle_i \qquad \text{if } n > i, d > 0,$$

(2.3)

$$\left\langle {n \atop d} \right\rangle_n = \left\langle {n-1 \atop d-1} \right\rangle \qquad \text{if } n = i, d > 0.$$

If $n = i = 1$, $\langle {n \atop d} \rangle_i$ is equal to 1 if $d = 0$ and equal to 0 if $d > 0$. The modified Eulerian number $\langle {n \atop d} \rangle_{i,j}$ is defined as the number of $\pi \in S_n$ with $\pi(1) = i$, $\pi(n) = j$ and $d$ descents. If $d = 0$, $\langle {n \atop d} \rangle_{i,j}$ is 1 if $i = 1$ and $j = n$ but 0 otherwise. For $d > 0$, the following recurrence can be derived:

$$\left\langle {n \atop d} \right\rangle_{i,j} = d \left\langle {n-1 \atop d} \right\rangle_{i,j} + (n-d-1) \left\langle {n-1 \atop d-1} \right\rangle_{i,j} \qquad \text{if } n > i, n > j,$$

(2.4)

$$= \left\langle {n-1 \atop d-1} \right\rangle_{n-j} \qquad \text{if } n = i, n > j,$$

$$= \left\langle {n-1 \atop d} \right\rangle_i \qquad \text{if } n > i, n = j.$$

If $n = i = j = 1$ and $d > 0$, $\langle {n \atop d} \rangle_{i,j} = 0$. Using (2.2), (2.3) and (2.4), the polynomials $p_{i,j}(x)$ can be formed in the base case.

The descent polynomial of shuffles from $1^{n_1}, 2^{n_2}, \ldots, h^{n_h}$ to $D$ is obtained as the sum of the polynomials $p_{i,j}(x)$ over $1 \leq i \leq n_1$ and $1 \leq j \leq n_h$.

We now turn to the descent polynomial of shuffles $\pi$ from $D$ to the sorted deck $1^{n_1}, 2^{n_2}, \ldots, h^{n_h}$. We first consider the occurrence of descents between $\pi(k)$ and $\pi(k+1)$ when the positions $k$ and $k+1$ are occupied in $D$ by cards with different labels. There will be a descent if and only if the label of the card at $k$ is greater than the label of the card at $k+1$. Thus, the number of descents of this type is the same for every shuffle from $D$ to the sorted deck and is equal to the number of places where a card with a greater label immediately precedes a card with a lesser label in the deck $D$. This quantity, which may be denoted by $\text{des}(D)$, is called the number of descents in $D$ and is extensively studied in [10] and [12].

We next consider descents between $\pi(k)$ and $\pi(k+1)$ only if both positions $k$ and $k+1$ are occupied in $D$ by cards with the label $c$. The cards at $k$ and $k+1$ both have the label $c$ if and only if $k = D(i,c)$ and $k+1 = D(i+1,c)$ for some integer $i$, $1 \leq i < n_c$. To facilitate the counting of this type of descent, denote the generating polynomial $\sum_{d=0}^{n-1} \langle {n \atop d} \rangle x^d$ of the Eulerian numbers by $\eta_n(x)$ [10]. If we pay attention only to cards with label $c$ in the deck $D$, it will look like $**ccc***cc**c$ with blocks of $c$'s separated by cards with different labels. Assume that the lengths of these blocks are given by $m_1, m_2, \ldots, m_\gamma$, with $\gamma$ being the number of blocks. Then $n_c = m_1 + m_2 + \cdots + m_\gamma$. Let

$$j_c = \sum_{j=1}^{c} n_j$$

and let $i_c = j_c - n_c + 1$. If $\pi$ is a shuffle from $D$ to the sorted deck, then $i_c \leq \pi(D(i,c)) \leq j_c$ must hold for $1 \leq i \leq n_c$. The $n_c$ integers in $[i_c, j_c]$ can be divided into sets of $m_1, m_2, \ldots, m_\gamma$ in $n_c!/(m_1! m_2! \cdots m_\gamma!)$ ways. For each such division of these $n_c$ integers into sets, there are $m_1! m_2! \cdots m_\gamma!$ ways of assigning values to $\pi(D(i,c))$, $1 \leq i \leq n_c$, such that a number assigned to a position within the first block of $c$s is in the first set and so on. The coefficient of $x^d$ of the polynomial $\eta_{m_1}(x) \eta_{m_2}(x) \cdots \eta_{m_\gamma}(x)$ is equal to the number of these assignments which have $d$ descents. Therefore, the coefficient of $x^d$ of the polynomial

$$(2.5) \qquad p_c(x) = \frac{n_c!}{m_1! m_2! \cdots m_\gamma!} \eta_{m_1}(x) \eta_{m_2}(x) \cdots \eta_{m_\gamma}(x)$$

is equal to the number of assignments with $d$ descents out of the $n_c!$ assignments of integers in $[i_c, j_c]$ to $\pi(D(i,c))$, $1 \leq i \leq n_c$. As intended, (2.5) counts the descent between $\pi(k)$ and $\pi(k+1)$ if and only if cards at positions $k$ and $k+1$ in $D$ both have the label $c$.

To find the descent polynomial of shuffles from $D$ to the sorted deck, note that the occurrence of a descent between $\pi(k)$ and $\pi(k+1)$, with cards labeled $c$ at positions $k$ and $k+1$ in $D$, is completely independent of the occurrence of a descent between $\pi(l)$ and $\pi(l+1)$, with cards labeled $d$ at positions $l$ and $l+1$, if $c \neq d$. Moreover, there are always $\text{des}(D)$ descents in a shuffle $\pi$ from $D$ to the sorted deck that correspond to positions $k$ and $k+1$ occupied in $D$ by cards with different labels. Therefore, the descent polynomial is given by

$$(2.6) \qquad x^{\text{des}(D)} p_1(x) p_2(x) \cdots p_h(x),$$

where the $p_i(x)$ are defined by (2.5).

If the deck $D$ is any permutation of the multiset $\{1^{n_1}, 2^{n_2}, \ldots, h^{n_h}\}$, (2.1) and (2.6) make it possible to find the descent polynomials of shuffles from the sorted deck to $D$ and of shuffles from $D$ to the sorted deck in polynomial time. The descent polynomial of shuffles between decks neither of which is sorted will be considered in later work.

In the rest of this section, we turn to theorems about transition probabilities between decks under an $a$-shuffle which do not use the descent polynomial. Let $a_1, a_2, \ldots, a_n$ be one of the $a^n$ integer sequences with $0 \leq a_i < a$ for $1 \leq i \leq n$. This sequence can be sorted to $a_{i_1} \leq a_{i_2} \leq \cdots \leq a_{i_n}$ in a stable manner and the permutation $i_1, i_2, \ldots, i_n$ of $\{1, 2, \ldots, n\}$ is uniquely determined since we require $i_j < i_{j+1}$ if $a_{i_j} = a_{i_{j+1}}$. Associate the shuffle $\pi \in S_n$ with $\pi(k) = i_k$ for $1 \leq k \leq n$ with the sequence $a_1, a_2, \ldots, a_n$. Then the uniform distribution on the $a^n$ sequences induces the $a$-shuffle distribution on $S_n$ [2]. This description of the $a$-shuffle is used in Theorems 2.1 and 2.2.

THEOREM 2.1. *Among all decks $D$ that are permutations of the multiset $\{1^{n_1}, 2^{n_2}, \ldots, h^{n_h}\}$, the transition probability under an $a$-shuffle from the sorted deck $1^{n_1}, 2^{n_2}, \ldots, h^{n_h}$ to $D$ is greatest for $D = 1^{n_1}, 2^{n_2}, \ldots, h^{n_h}$ and least for $D = h^{n_h}, \ldots, 2^{n_2}, 1^{n_1}$.*

PROOF. Assume the sorted deck to be $1^n, 2^n$. The proof for more general sorted decks is similar.

Let $a_1, a_2, \ldots, a_{2n}$ be a sequence with $0 \leq a_i < a$ for $1 \leq i \leq 2n$. If $D(i, 1) = j$, define $\alpha_i = a_j$, and if $D(i, 2) = j$, define $\beta_i = a_j$. For example, if $D = 1, 2, 1, 2, 1, 2$, then

$$a_1, a_2, a_3, a_4, a_5, a_6 = \alpha_1, \beta_1, \alpha_2, \beta_2, \alpha_3, \beta_3.$$

For the sequence to induce a shuffle from $1^n, 2^n$ to $D$, each $\alpha_i$ must be less than or equal to each $\beta_i$. In addition, each $\alpha_i$ must be strictly less than all the $\beta$'s that precede it in the sequence. For example, if $D = 1, 2, 1, 2, 1, 2$, the inequalities are

$$\max(\alpha_1, \alpha_2, \alpha_3) \leq \min(\beta_1, \beta_2, \beta_3), \qquad \alpha_2 < \beta_1, \qquad \alpha_3 < \beta_1, \qquad \alpha_3 < \beta_2.$$

If $D = 1^n, 2^n$, it is enough if each $\alpha$ is less than or equal to each $\beta$. If $D = 2^n, 1^n$, each $\alpha$ must be strictly less than each $\beta$. Therefore, the number of sequences that induce a shuffle from $1^n, 2^n$ to $D$ is greatest for $D = 1^n, 2^n$ and least for $D = 2^n, 1^n$. The statement about transition probabilities follows. □

Theorem 2.2 below generalizes Theorem 3 of [2] and their proofs use similar arguments. Similar arguments can also be found in [6] and [10].

THEOREM 2.2. *Let the deck $D$ be a permutation of the multiset $\{1, 2, \ldots, h, x^n\}$. Let the number of cards labeled $c$, $1 \leq c \leq h$, that are not preceded by a card labeled $c - 1$ in $D$ be equal to $r$. Let the number of cards labeled $x$ that precede the card labeled $h$ in $D$ be equal to $l$. Then the probability that an $a$-shuffle applied to the sorted deck $1, 2, \ldots, h, x^n$ results in $D$ is*

$$\frac{1}{a^{n+h}} \sum_{m=r-1}^{a-1} \binom{m-r+h}{h-1} (a-m-1)^l (a-m)^{n-l},$$

*where if $l = 0$ and $m = a - 1$, $(a - m - 1)^l$ must be taken to be 1.*

PROOF. Let $a_1, a_2, \ldots, a_{h+n}$ be an integer sequence with $0 \leq a_i < a$ for $1 \leq i \leq h + n$. If $D(1, c) = i$, $1 \leq c \leq h$, define $\alpha_c = a_i$. If $D(j, x) = k$, $1 \leq j \leq n$, define $\beta_j = a_k$. For the sequence $a_1, a_2, \ldots, a_{h+n}$ to induce an $a$-shuffle from the sorted deck to $D$, we require

$$(2.7) \qquad \alpha_1 \leq \alpha_2 \leq \cdots \leq \alpha_h \leq \min(\beta_1, \ldots, \beta_n).$$

In addition, the inequality $\alpha_{c-1} \leq \alpha_c$, $2 \leq c \leq h$, must be strict if the card labeled $c$ in $D$ is not preceded by the card labeled $c-1$. Therefore, exactly $r-1$ inequalities between the $\alpha$s in (2.7) are strict. Further, at least $l$ of the $\beta_i$, the ones with $1 \leq i \leq l$, must be strictly greater than $\alpha_h$.

The number of solutions to (2.7), with the additional conditions described below it, can be counted by allowing $\alpha_h = m$ to vary from $r-1$ to $a-1$. Given $m$, the number of ways to pick the $\alpha_c$, $1 \leq c < h$, can be counted as follows. Start with $m$ "jumps." Allocate $r-1$ of these jumps to the inequalities in $\alpha_1 \leq \alpha_2 \leq \cdots \leq \alpha_{h-1} \leq m$ that must be strict. The remaining $m - r + 1$ jumps can be assigned to $h$ positions, namely, the position before $\alpha_1$ and the $h-1$ inequalities, in $\binom{m-r+h}{h-1}$ ways. The value of $\alpha_i$, for $1 \leq i \leq h-1$, is equal to the number of jumps preceding it. The number of ways to pick the $\beta$s is $(a - m - 1)^l (a - m)^{n-l}$. The formula for the transition probability from the sorted deck to $D$ follows. $\square$

**3. Starting deck $1^n, 2^n$.** The probability distribution over decks that are permutations of the same multiset of cards under an $a$-shuffle can be obtained from (2.1) or (2.6) if either the starting deck or the ending deck is in sorted order. The total variation distance from the uniform distribution is a sum over all possible decks and its calculation can therefore involve a very large number of terms. However, the calculation becomes simpler if it is recognized that the transition probabilities are the same for whole classes of decks. In the case where all $n$ cards have different labels, the transition probabilities depend only upon the number of descents in the shuffle and, hence, the $n!$ decks fall into only $n$ equivalence classes. In this section we investigate this type of equivalence relationship when the starting deck is $1^n, 2^n$.

In this section and the next, we use $\alpha$, $\beta$ and $\gamma$ to denote sequences of 1's and 2's that stand for segments of a deck of cards. The number of entries in the sequence $\alpha$ is denoted by $|\alpha|$. The sequence obtained by reversing the order of $\alpha$ and then replacing each 1 by 2 and each 2 by 1 is denoted $\alpha^*$. For example, if $\alpha = 1, 2, 2, 2, 1, 1$, then $\alpha^* = 2, 2, 1, 1, 1, 2$. A total of $\binom{2n}{n}$ decks can be obtained by rearranging the cards of $1^n, 2^n$. The equivalence relation $R$ on that set of decks is defined as follows. The deck $D_1 = \alpha\beta\gamma$ is $R$-related to $D_2 = \alpha\beta^*\gamma$ if $|\alpha| = |\gamma|$, and the number of 1's and the number of 2's in $\beta$ are equal. For example, $1, 2, 2, 1$ is $R$-related to $2, 1, 1, 2$ and $1, 1, 2, 2, 1, 2$

is $R$-related to $1, 2, 1, 1, 2, 2$. The equivalence relation is obtained by taking the transitive, reflexive closure. For example, the decks $1, 2, 1, 2, 2, 1, 2, 1$ and $1, 2, 2, 1, 1, 2, 2, 1$ and $2, 1, 1, 2, 2, 1, 1, 2$ and $2, 1, 2, 1, 1, 2, 1, 2$ are all in the same equivalence class.

THEOREM 3.1. *If $D_1$ is $R$-related to $D_2$, the transition probability from $1^n, 2^n$ to $D_1$ is equal to the transition probability from $1^n, 2^n$ to $D_2$ under an $a$-shuffle for any $a$.*

PROOF. It is sufficient to consider $D_1 = \alpha\beta\gamma$ and $D_2 = \alpha\beta^*\gamma$ with $|\alpha| = |\gamma|$ and with equal number of 1's and 2's in $\beta$. It is enough to show that the descent polynomial of shuffles from $1^n, 2^n$ to $D_1$ is equal to the descent polynomial of shuffles from $1^n, 2^n$ to $D_2$. We will construct a bijective map from $\pi(1^n, 2^n; D_1)$ to $\pi(1^n, 2^n; D_2)$ such that a shuffle maps to another shuffle with exactly the same number of descents.

The number of 1's in $\alpha$ is equal to the number of 2's in $\gamma$ since the number of 1's and 2's are equal in $\alpha\beta\gamma$ and in $\beta$. Similarly, the number of 2's in $\alpha$ is equal to the number of 1's in $\gamma$. Let $p_1 < \cdots < p_a$ and $q_1 < \cdots < q_b$ and $r_1 < \cdots < r_c$ be the positions of $D_1$ that correspond to $\alpha$ and $\beta$ and $\gamma$, respectively, that are occupied by 1's. Similarly, let $p'_1 < \cdots < p'_c$ and $q'_1 < \cdots < q'_b$ and $r'_1 < \cdots < r'_a$ be the positions of $D_1$ that correspond to $\alpha$ and $\beta$ and $\gamma$, respectively, that are occupied by 2's.

Define $f(x) = 2n - x + 1$. The map $f$ reflects a position in a deck of size $2n$ about its center; for example, the first position is reflected to the last position. The positions occupied by 1's in $D_2 = \alpha\beta^*\gamma$ are

$$p_1 < \cdots < p_a < f(q'_b) < \cdots < f(q'_1) < r_1 < \cdots < r_c,$$

where $p$s and $r$s correspond to $\alpha$ and $\gamma$, and $f(q)$s correspond to $\beta^*$. In $D_1$, the position $q'_i$ is occupied by a 2. When $\beta$ is reversed that 2 is moved to the position $f(q'_i)$ and then it is replaced by 1 to form $\beta^*$. This explains the central block of $f(q'_i)$'s above. Similarly, the positions occupied by 2's in $D_2$ are

$$p'_1 < \cdots < p'_c < f(q_b) < \cdots < f(q_1) < r'_1 < \cdots < r'_a,$$

where the $p'$'s correspond to positions in $\alpha$, $f(q)$'s to positions in $\beta^*$ and $r'$'s to positions in $\gamma$. Note that each $p$ or $p'$ is less than each $q$ or $q'$, which is less than each $r$ or $r'$.

Let $\pi \in S_{2n}$ be a shuffle from $1^n, 2^n$ to $D_1$. Then the numbers

$$\pi(1), \pi(2), \ldots, \pi(n)$$

must be an arrangement of the positions occupied by 1's in $D_1$. Similarly, the numbers

$$\pi(n+1), \pi(n+2), \ldots, \pi(2n)$$

must be an arrangement of the positions occupied by 2's in $D_1$. The map to a shuffle from $1^n, 2^n$ to $D_2$ is based on two cases. In the first case, we assume that *not both* $\pi(n)$ and $\pi(n+1)$ correspond to positions in $\beta$. The shuffle $\pi^*$ from $1^n, 2^n$ to $D_2$ that $\pi$ maps to is defined as

$$\pi^*(i) = \phi(\pi(i)),$$

for $1 \leq i \leq 2n$, where $\phi(\cdot)$ will now be defined. First, we define $\phi(p_i) = p_i$, $\phi(r_i) = r_i$, $\phi(p'_i) = p'_i$ and $\phi(r'_i) = r'_i$. In addition, we define $\phi$ as

$$q_1 \to f(q'_b) \quad q_2 \to f(q'_{b-1}) \quad \cdots \quad q_b \to f(q'_1)$$
$$q'_1 \to f(q_b) \quad q'_2 \to f(q_{b-1}) \quad \cdots \quad q'_b \to f(q_1).$$

This definition maps the $q$s to the $f(q')$'s and the $q'$'s to the $f(q)$'s and, therefore, $\pi^*$ is a shuffle from $1^n, 2^n$ to $D_2$. Further, $x < y$ if and only if $\phi(x) < \phi(y)$, except when $x$ is a $q$ and $y$ is a $q'$ or when $x$ is a $q'$ and $y$ is a $q$. However, in the arrangement

$$\pi(1), \ldots, \pi(n), \pi(n+1), \ldots, \pi(2n),$$

a $q$ and $q'$ can occur in consecutive positions as $\pi(n)$ and $\pi(n+1)$, and in no other way, and we have assumed that *not both* of those positions correspond to $\beta$. Therefore, the above arrangement has the same number of descents as

$$\phi(\pi(1)), \ldots, \phi(\pi(n)), \phi(\pi(n+1)), \ldots, \phi(\pi(2n))$$

and $\pi^*$ has the same number of descents as $\pi$.

The other case is when $\pi(n)$ is a $q$ and $\pi(n+1)$ is a $q'$. Then we define $\phi(q_i) = f(q_i)$, $\phi(q'_i) = f(q'_i)$, $p_i \xleftrightarrow{\phi} r'_{a-i+1}$ and $r_i \xleftrightarrow{\phi} p'_{c-i+1}$. We map $\pi$ to $\pi^*$, where $\pi^*$ is defined as

$$\pi^*(i) = \phi(\pi(2n - i + 1)),$$

for $1 \leq i \leq 2n$. It can be verified that $\pi^*$ is a shuffle from $1^n, 2^n$ to $D_2$. Further, $x < y$ if and only if $\phi(x) > \phi(y)$, except when $x$ is a $p$ and $y$ is a $p'$, or $x$ is a $p'$ and $y$ is a $p$, or $x$ is an $r$ and $y$ is an $r'$, or $x$ is an $r'$ and $y$ is an $r$. In the arrangement,

$$\pi(1), \ldots, \pi(n), \pi(n+1), \ldots, \pi(2n),$$

a $p$ and $p'$ or an $r$ and $r'$ can occur in consecutive positions only at $\pi(n)$ and $\pi(n+1)$. However, we have assumed that $\pi(n)$ is a $q$ and that $\pi(n+1)$ is a $q'$. Therefore, every descent in the above arrangement becomes an ascent in

$$\phi(\pi(1)), \ldots, \phi(\pi(n)), \phi(\pi(n+1)), \ldots, \phi(\pi(2n))$$

and every ascent becomes a descent. This arrangement is reversed to define $\pi^*(1), \ldots, \pi^*(2n)$ which changes the ascents back into descents and, therefore, the number of descents in $\pi^*$ is equal to the number of descents in $\pi$.

Finally, we need to show that the map defined above is a bijection. A shuffle from $1^n, 2^n$ to $D_2$ can be mapped to a shuffle from $1^n, 2^n$ to $D_1$ using the same procedure as above. The resulting map is the inverse of the above map because $\phi \circ \phi$ is identity in both cases above. $\square$

It is natural to ask if the equality of the descent polynomials of the shuffles from $1^n, 2^n$ to $D_1$ and $D_2$ implies that $D_1$ is $R$-related to $D_2$. We have checked that this is indeed so for $n = 1, 2, 3, 4, 5, 6, 7$. The theorem below counts the total number of equivalence classes under the relation $R$.

THEOREM 3.2. *The number of equivalence classes under $R$ is equal to the Catalan number $\frac{1}{n+2}\binom{2n+2}{n+1}$.*

PROOF. We describe a method to find a unique representative for each equivalence class and then count the number of unique representatives. The function $f(x) = 2n - x + 1$ reflects positions with respect to the center of the deck as before. In this proof, we refer to $f(x)$ as the reflection of the position $x$. A position and its reflection either both lie in $\beta$ or lie outside it, since $|\alpha| = |\gamma|$. Consider the positions $x = n+1, n+2, \ldots, 2n$. Assume that $x$ and $f(x)$ both lie inside $\beta$. If the positions $x, f(x)$ are occupied by $1, 2$ in the deck $D$, reversal of $\beta$ changes it to $2, 1$, and the replacements of 1's by 2's and 2's by 1's changes it back to $1, 2$. Similarly, application of the basic rule that generates the relation $R$ does not change $D$ in the positions $x, f(x)$ if those positions are occupied by $2, 1$. However, if they are occupied by $1, 1$, that becomes $2, 2$ when the rule is applied using a $\beta$ large enough to include positions $x$ and $f(x)$. Similarly, $2, 2$ becomes $1, 1$. With each position $x$, we associate the symbol "+" if positions $x, f(x)$ are occupied by $1, 2$, the symbol "−" if occupied by $2, 1$, the symbol 1 if occupied by $1, 1$, and the symbol 2 if occupied by $2, 2$. The deck as a whole is coded as the list of symbols associated with positions $n+1$ through $2n$. For example, the deck $1, 1, 2, 1, 2, 2, 2, 1, 1, 2, 1, 2, 1, 2$ is coded as $+, +, 2, 1, 2, 1, -$.

The +s and −s never change when the rule that generates the relation $R$ is repeatedly applied with possibly many different choices of $\beta$. They are ignored in much of the rest of this proof. We can find the $\beta$'s which lead to a nontrivial application of the rule to generate the relation $R$ using the code for the deck $D$ as follows. We traverse the code from left to right, and record the excess of 1's over 2's. For example, for the code $-, -, 2, 1, 2, 1, +$, this excess is $-1$ after the first 2 is passed, then becomes 0, and then $-1$, and then 0. The rule for generating $R$ can be applied whenever this excess becomes 0. If the excess becomes zero, after traversing $i$ symbols in the code, the corresponding $\beta$ in the deck is a segment of $2i$ cards extending from position $n-i+1$ to position $n+i$. When the rule is applied, the 1's become 2's and the 2's become 1's among the first $i$ symbols of the code. If the excess becomes 0 after $i$ symbols and then again after $i+j$ symbols, the application of the rule with a $\beta$ of length equal to $2i$ followed by another application using a $\beta$ of length $2(i+j)$ changes the code for the deck only between the $(i+1)$st and the $j$th symbol. Among these symbols, the 1's change to 2's and the 2's change to 1's. By applying the rule with judicious choices of $\beta$, it is possible to obtain a single code in which the excess never becomes negative. For example, the code $-,-,2,1,2,1,+$ can be converted to $-,-,1,2,1,2,+$. We use such codes as unique representatives of equivalence classes of decks.

Assume that in such a code, there are $k$ symbols equal to 1 and $k$ symbols equal to 2. Then there must be $n-2k$ symbols equal to $+$ or $-$. If the $+$ and $-$ symbols are ignored, and each 1 is substituted by a ( and each 2 by a ), we obtain a valid arrangement of parentheses of length $2k$. The number of valid arrangements of parentheses of length $2k$ is well known to be the Catalan number $\frac{1}{k+1}\binom{2k}{k}$. For each assignment of 1's and 2's to the $2k$ positions, the other positions can be filled with symbols $+$ and $-$ in $2^{n-2k}$ ways. The $2k$ positions that are assigned either the symbol 1 or the symbol 2 can be chosen in $\binom{n}{2k}$ ways. Therefore, the total number of equivalence classes is given by $\sum_k \frac{2^{n-2k}}{k+1}\binom{2k}{k}\binom{n}{2k}$. This sum can be simplified:

$$\sum_k \frac{2^{n-2k}}{k+1}\binom{2k}{k}\binom{n}{2k} = 2^n \sum_k \frac{1}{k+1}\binom{n/2}{k}\binom{(n-1)/2}{k}$$
$$= \frac{2^n}{n/2+1}\sum_k \binom{n/2+1}{k+1}\binom{(n-1)/2}{k}.$$

The first equality above uses (5.35) in [9]. The proof may be completed using the binomial identity $\sum_k \binom{l}{m+k}\binom{s}{p+k} = \binom{l+s}{l-m+p}$ for integers $l,m,p$ and $l \geq 0$. The cases with $n$ even and odd have to be considered separately. □

**4. Starting deck $(1,2)^n$.** The equivalence relation $R$ in this section is different from the one considered in the previous section. In this section $D_1 = \alpha\beta\gamma$ is $R$-related to $D_2 = \alpha\beta^*\gamma$ if $\beta$ has the same number of 1's as 2's. The additional condition $|\alpha| = |\gamma|$ is no longer required. The decks $D_i$ are all permutations of $\{1^n, 2^n\}$. The equivalence relation is obtained by taking the transitive, reflexive closure. The equivalence class containing $1,2,1,2,2,1,2,1$ has five other decks in it.

THEOREM 4.1. *If $D_1$ is $R$-related to $D_2$, the transition probability from $(1,2)^n$ to $D_1$ is equal to the transition probability from $(1,2)^n$ to $D_2$ under an $a$-shuffle for any $a$.*

PROOF. It is sufficient to consider $D_1 = \alpha\beta\gamma$ and $D_2 = \alpha\beta^*\gamma$ with equal number of 1's and 2's in $\beta$. We will construct a bijective map from $\pi((1,2)^n; D_1)$ to $\pi((1,2)^n; D_2)$ such that a shuffle maps to another shuffle with exactly the same number of descents.

Let $\pi \in \pi((1,2)^n; D_1)$. Let $\pi(2i-1) = a_i$, $1 \leq i \leq n$, and $\pi(2i) = b_i$, $1 \leq i \leq n$. The number of descents in $\pi$ is equal to the number of descents in the arrangement $a_1, b_1, a_2, b_2, \ldots, a_n, b_n$. To facilitate the proof, we depict this arrangement in the following way:

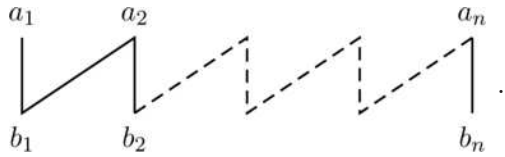

In the deck $D_1$ each position $a_i$ is occupied by a 1 and each position $b_i$ is occupied by a 2 [because $\pi$ is a shuffle from $(1,2)^n$ to $D_1$]. We assume that $\beta$ begins at the $(i+1)$st position and ends at the $(i+2j)$th position. In the depiction above, circle all the $a_i$ and $b_i$ that do not correspond to $\beta$, that is, circle an $a_i$ or a $b_i$ if it is less than $i+1$ or greater than $i+2j$. When $\pi$ is mapped, the circled numbers will stay fixed. The uncircled numbers form segments that run from a circled number to another circled number (or they might begin or end at $a_1$ or $b_n$). These segments of uncircled numbers are of four types according as they either begin or end at the top line or the bottom line:

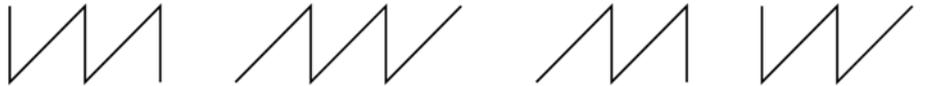

The third type of segment has one more uncircled position in the bottom line, corresponding to a position occupied by a 2 in $D_1$, than in the top line. The fourth type of segment has an extra uncircled position in the top line, corresponding to a position occupied by a 1 in $D_1$. Since the number of 1's in $\beta$ is equal to the number of 2's, the number of uncircled positions in the top line must be equal to the number of uncircled positions in the bottom line. Therefore, the number of uncircled segments of the third type must be equal to the number of uncircled segments of the fourth type.

To map $\pi$ to a shuffle from $(1,2)^n$ to $D_2$, we will modify the uncircled segments and insert them back in-between the circled numbers in the original arrangement of $a_i$ and $b_i$. We define a map $f$ from the uncircled positions, that is, the positions that correspond to $\beta$, back to the the uncircled positions as follows:

$$i+1 \to i+2j,\ i+2 \to i+2j-1,\ \ldots,\ i+2j-1 \to i+2,\ i+2j \to i+1.$$

If $i+1 \leq x \leq i+2j$ and the position $x$ in $D_1$ is occupied by 1 (or 2), the position $f(x)$ in $D_2$ will be occupied by 2 (or 1). If $a_p, b_p, a_{p+1}, b_{p+1}, \ldots, b_q$ is an uncircled segment of the first type, it will be modified to $f(b_q), f(a_q), \ldots, f(a_{p+1}), f(b_p), f(a_p)$. In the deck $D_2$, each position $f(b_i)$ [or $f(a_i)$], $p \leq i \leq q$, is occupied by 1 (or by 2). Therefore, the modified segment is also of the first type. However, when an uncircled segment of the third (or fourth) type is modified in this way, it becomes a segment of the fourth (or third) type. The arrangement $a_1, b_1, a_2, b_2, \ldots, a_n, b_n$ can be converted to another arrangement in the following steps:

1. Extract the uncircled segments of the first and second type from the arrangement, modify them as described above, and put the modified segment back in the same place.
2. Number the uncircled segments of the third and fourth type from left to right. As explained above, they must be equally numerous.
3. Replace the $i$th uncircled segment of the third type by the modification of the $i$th uncircled segment of the fourth type. Similarly, replace the $i$th uncircled segment of the fourth type by the modification of the $i$th uncircled segment of the third type.

The shuffle $\pi^*$ is defined by setting $\pi^*(i)$ equal to the $i$th number in the arrangement constructed in this manner. By construction, $\pi^*$ is a shuffle from $(1,2)^n$ to $D_2$. Further, the number of descents of $\pi^*$ must equal the number of descents of $\pi$ for the following reason. If there is a descent or an ascent between two circled positions in the arrangement $a_1, b_1, a_2, b_2, \ldots, a_n, b_n$, it remains unchanged. Also, the modification of uncircled segments described above preserves the number of descents, although it changes their locations. Finally, if a circled number is greater than (or less than) a single uncircled number, it must be greater than (or less than) all uncircled numbers and, therefore, the number of descents between circled and uncircled numbers in the arrangement also remains unchanged.

It is possible to map a shuffle from $(1,2)^n$ to $D_2$ to a shuffle from $(1,2)^n$ to $D_1$ using the same procedure. That map would be the inverse of the map defined above. Therefore, the map from shuffles $\pi$ to shuffles $\pi^*$ defined above is a bijection.  $\square$

The converse of the above Theorem 4.1 appears to be true as well. The number of equivalence classes seems to be given by the simple formula $(n+3)2^{n-2}$. One can attempt to prove this by finding unique representatives for equivalence classes and then counting them as in the proof of Theorem 3.2. We have derived a method to construct unique representatives for equivalence classes of $R$. However, we have not yet devised a method to count the number of unique representatives.

**5. Card games.** Some inferences about the mixing times for common card games such as blackjack and bridge can be drawn using results given in the preceding sections. Let $S$ be a finite set and let $p$ be a probability distribution on $S$. Then the total variation distance of $p$ from the uniform distribution is given by $\frac{1}{2}\sum_{s\in S}|p(s) - \frac{1}{|S|}|$. For a deck of 52 distinct cards, the total variation distance remains close to 1 until the number of riffle shuffles exceeds 4. The total variation distance falls below 0.5 when the number of shuffles is 7 and this can be taken to be the mixing time [2]. Another distance defined in [1] is the separation distance. The separation distance of $p$ from the uniform distribution is $\max_{s\in S}(1 - |S|p(s))$. Like the total variation distance, the separation distance has a maximum of 1 and a minimum of 0. However, it leads to a more demanding notion of mixing as the number of riffle shuffles of a deck of 52 distinct cards needed to make the separation distance no more than $1/2$ is 11. The use of entropy to understand mixing is discussed in [14]. The validity and limitations of the maximum entropy model of riffle shuffles are discussed in [3] and [5].

In the game of bridge, 52 distinct cards are dealt to four players. To apply the results proved in Section 2, we need to assume that the first 13 cards are dealt to one player, the next 13 to another and so on. Let the deck $D$ be a permutation of the multiset $\{1^{13}, 2^{13}, 3^{13}, 4^{13}\}$. Let $p_D$ be the transition probability from $D$ to $1^{13}, 2^{13}, 3^{13}, 4^{13}$ under an $a$-shuffle. This transition probability can be obtained using (2.6). The probability that the first player is dealt cards originally in the positions occupied by cards labeled 1 in $D$, that the second player is dealt cards originally in the positions occupied by cards labeled 2 in $D$, and so on after an $a$-shuffle is equal to $p_D$. Therefore, the distance of the probability distribution $p_D$ over decks from the uniform distribution will indicate the closeness of a deal after an $a$-shuffle to a random deal to four players. If the separation distance is used to define mixing, an application of (2.6) with $D = (4, 3, 2, 1)^{12}$ shows that the separation distance is greater than 0.5 after 10 riffle shuffles. Therefore, the mixing time is still 11 riffle shuffles. The total variation distance involves a sum with a great number of terms and the results of Sections 3 and 4 indicate that the recognition of equalities among the transition probabilities $p_D$ is unlikely to make this sum tractable. However, a Monte Carlo procedure for evaluating this sum, which will be described elsewhere, implies that the mixing time is 6 riffle shuffles when the total variation distance is used. If the cards are dealt to the players in cyclic order, which is the common practice, the mixing times will almost certainly be lower.

In the game of blackjack, the distinction between the suits is ignored. We assume the starting deck to be $1^4, 2^4, \ldots, 13^4$. Application of Theorem 2.1 and (2.1) shows that the separation distance from the uniform distribution over decks becomes less than 0.5 after 9 riffle shuffles. Again, a Monte-Carlo procedure has to be employed to find the total variation distance. It then follows that the total variation distance becomes less than 0.5 after only 4 riffle shuffles.

**Acknowledgments.** The authors thank Professors P. Diaconis, J. Fulman, S. Lalley, C. Mulcahy and J. Rauch for helpful discussions.

DEPARTMENT OF MATHEMATICS
UNIVERSITY OF MICHIGAN
525 E. UNIVERSITY AVENUE
ANN ARBOR, MICHIGAN 48109
USA
E-MAIL: mconger@umich.edu
         divakar@umich.edu